\documentclass{article}
\usepackage{graphicx}
\usepackage{booktabs}
\usepackage[a4paper, total={5.5in, 9in}]{geometry}
\usepackage{array}
\usepackage{caption}
\usepackage{float}
\usepackage{multirow}
\usepackage{amsmath}
\usepackage{amssymb}
\usepackage{setspace}
\usepackage{pifont}
\usepackage{breqn}
\usepackage{multicol}
\usepackage{color}
\usepackage[english]{babel}
\usepackage[round]{natbib}
\usepackage{hyperref}

\usepackage{geometry} 
\usepackage{algorithm}
\usepackage{algorithmic}
\usepackage{tikz}
\usepackage{pgfgantt}
\usepackage{wrapfig}
\usepackage{adjustbox}
\usepackage{amsthm}
\usepackage{pst-node}
\usepackage{pstricks}
\usepackage{colortbl}
\usepackage{tabulary}
\usepackage{subcaption}
\usepackage{algorithmic}
\usepackage{lscape}
\usepackage{float}
\usepackage{booktabs}
\usepackage{textgreek}
\usepackage{threeparttable}
\usepackage{pdfpages}
\usepackage{caption}
\usepackage[title,toc,titletoc,page]{appendix}

\newcommand{\ignore}[1]{}

\providecommand{\keywords}[1]{\textbf{\textit{Keywords:}} #1}

\theoremstyle{definition}
\newtheorem{definition}{Definition}[section]

\title{{\Huge Bi-objective facility location in the presence of uncertainty}}
\author{\Large Najmesadat Nazemi $^a$, Sophie N. Parragh $^a$, Walter J. Gutjahr $^b$ \vspace{3ex} \\ $^a$ Institute of Production and Logistics Management, Johannes Kepler University Linz\\Altenberger Straße 69, 4040 Linz, Austria\\\href{najmesadat.nazemi@jku.at}{\{najmesadat.nazemi, }\href{sophie.parragh@jku.at}{sophie.parragh\}@jku.at} \vspace{3ex}\\ 
$^b$ Department of Statistics and Operations Research, University of Vienna\\
	Oskar-Morgenstern-Platz 1, 1090 Vienna, Austria\\
	\href{walter.gutjahr@univie.ac.at}{walter.gutjahr@univie.ac.at}}
\date{}
\begin{document}
\maketitle

\vspace{10pt}
\begin{abstract}
	
	 Multiple and usually conflicting objectives subject to data uncertainty are main features in many real-world problems. Consequently, in practice, decision-makers need to understand the trade-off between the objectives, considering different levels of uncertainty in order to choose a suitable solution.
	 In this paper, we consider a two-stage bi-objective single source capacitated model \textcolor{black}{as a base formulation for designing} a last-mile network in disaster relief where one of the objectives is subject to demand uncertainty. We analyze scenario-based two-stage risk-neutral stochastic programming, adaptive (two-stage) robust optimization, and a two-stage risk-averse stochastic approach using conditional value-at-risk (CVaR). To cope with the bi-objective nature of the problem, we embed these concepts into two criterion space search frameworks, the $\epsilon$-constraint method and the balanced box method, to determine the Pareto frontier. Additionally, a matheuristic technique is developed to obtain high-quality approximations of the Pareto frontier for large-size instances. \textcolor{black}{In an extensive computational experiment}, we evaluate and compare the performance of the applied approaches based on real-world data from a Thies drought case, Senegal.\\
	\keywords{Multi-objective optimization, Uncertainty, Facility location problem, Humanitarian relief}
	
\end{abstract}

\pagenumbering{gobble}

\newpage
\pagenumbering{arabic}

\section{Introduction}
\label{sec:intro}

Decision-makers (DMs) often face multiple goals, which are in conflict with each other \citep{ehrgott2005multicriteria}, e.g., minimization of cost vs. maximization of service level, and, finding a feasible solution that simultaneously optimizes all criteria is usually impossible. Consequently, it is important for DMs to understand the trade-off between the considered objectives. To cope with this issue, many real-world problems have been modeled as multi-objective optimization (MOO) problems and, a variety of algorithms have been developed to produce the set of trade-off solutions.

Moreover, usually, DMs do not make their decisions in a completely certain environment. Depending on the problem, different levels of uncertainty should be taken into account to make reliable decisions. To deal with this issue in optimization problems, different approaches have been proposed in the literature. The most widely used ones are stochastic programming \citep{birge2011introduction} and robust optimization \citep{ben2009robust}. Stochastic optimization assumes that the underlying probability distributions of the uncertain parameters are known or can be estimated, e.g., based on historical data. It allows to model different risk attitudes, e.g., risk-neutral decision making or risk-averse decision making. Unlike stochastic optimization, robust optimization does not assume any probabilistic data but only uncertain parameters stemming from an uncertainty set (e.g., scenarios). Well-known robust concepts are minmax/static robustness and adaptive robust optimization \citep{bertsimas2011theory}. 

Although many different approaches have been developed to cope with multiple objectives and uncertainty separately in optimization problems, the intersection of these two main domains, i.e., the analysis of decision problems involving multiple objectives and parameter uncertainty simultaneously, only received comparably little attention. Combinations of MOO with stochastic programming concepts (\cite{abdelaziz2012solution, gutjahr2016stochastic}; \cite{charkhgard2020bi}) and robustness concepts \citep{ehrgott2014minmax,ide2016robustness} give a flavor of the diverse possible ways to cope with real-life decision making problems.

Optimization problems arising along the humanitarian relief chain (HRC) are real-world applications that motivated us in combining bi-objective optimization and optimization under uncertainty, as multiple objectives and inherent risk are distinguishing features of humanitarian logistics. HRC problems concern multi-stakeholder decision making in which a population of stakeholders seeks to balance their conflicting objectives and priorities. The objectives can be categorized into three main groups of criteria, efficiency criteria, effectiveness criteria, and equity criteria \citep{gralla2014assessing}. Besides, in a disaster situation, most of the information received at the disaster management center, such as the number of injured people, the amount of demand of the affected people, traveling times according to the network conditions and available commodities, etc. are inherently imprecise and uncertain. To avoid any inefficiency in the final selected solutions, such inherent uncertainty in input data should be taken into account. Thus, it is desirable to have optimization models and solution techniques to deal with the multi-objective nature and the uncertainty features of HRC simultaneously.

Disaster management operations are usually classified into four phases: mitigation and preparedness (pre-disaster), response, and recovery (post-disaster). These four phases together are known as the disaster management cycle \citep{altay2006or}. The preparedness and response are relief phases, whereas the mitigation and recovery are development phases. When a disaster (either a sudden-onset disaster that occurs as a single, distinct event, such as an earthquake or a slow-onset disaster that emerges gradually over time, such as a drought) strikes, relief goods are transported to the points of distributions (PODs) by relief organizations. Then, affected people walk or drive to these centers to collect their relief aid. Since the numbers and the locations of the PODs and their capacity, directly affect the performance of the relief chain in terms of response time and costs \citep{balcik2008facility}, facility location decisions in disaster management are of vital importance. This paper deals with a two-stage bi-objective \textcolor{black}{single source capacitated facility location (Bi-SSCFL)} problem where one objective is affected by data uncertainty.

We analyze traditional two-stage risk-neutral stochastic and adaptive (two-stage) robust optimization as well as a two-stage risk-averse stochastic framework which can be seen as a method between the other two approaches using a widely applied risk measure, namely the conditional value-at-risk (CVaR($\alpha$)). With the help of uncertainty level $\alpha$ ($\alpha \in [0,1)$), it can be adapted to the DM's risk preference. Two-stage stochastic and adaptive robust optimization are two special cases of the CVaR with $\alpha=0$ and $\alpha \approx 1$, respectively. Towards that end, we use the classical linear representation of CVaR. \textcolor{black}{Embedded into the epsilon-constraint framework and the balanced box method, we solve small and medium sized instances to optimality, and to solve larger instances,} we develop a matheuristic method which finds high-quality approximations of the set of trade-off solutions \textcolor{black}{and can be used to solve instances.}

The remainder of this paper is organized as follows. A literature review on uncertain multi-objective problems in the context of HRC is presented in Section~\ref{sec:lit}. In Section~\ref{sec:problemdec}, we define the uncertain bi-objective facility location problem addressed in this paper and its mathematical formulation. After that, in Section \ref{sec:method}, we summarize the developed solution approaches to solve the problem. Section \ref{result} presents and discusses the computational results. Finally, in Section \ref{sec:conclusion} we present conclusion and address potential future work.

\section{Related work}
\label{sec:lit}

As mentioned in the \nameref{sec:intro}, multiple objectives and uncertainty are two major characteristics of humanitarian decision support systems in practice. Uncertainty, as well as multi-objective optimization, have been considered abundantly in the literature on humanitarian relief (see \cite{hoyos2015or}, \cite{grass2016two}, \citet{gutjahr2016multicriteria}). However, let us mention that also in other application areas of Logistics and Supply Chain Management, the issues of uncertainty and multiple objectives play an important role. According to \citet{gutjahr2016multicriteria}, capturing uncertainty in multi-criteria optimization is critical in practice and still a comparably young field.

The majority of the studies in the literature that addresses MOO under uncertainty in different phases of disaster relief (e.g., \cite{tzeng2007multi}, \cite{zhan2011multi},\cite{tricoire2012bi}, \cite{najafi2013multi}, \cite{rezaei2014robust}, \cite{rath2016bi}, \cite{haghi2017developing}, \cite{ liu2017robust}, \cite{kinay2019multi}, \cite{parragh2020branch}) use either scenario-based traditional two-stage stochastic programming (e.g., \cite{zhan2011multi}, \cite{tricoire2012bi}, \cite{rath2016bi}, \cite{parragh2020branch}) or scenario-based worst-case robust optimization (e.g., \cite{najafi2013multi}, \cite{rezaei2014robust}, \cite{haghi2017developing}, \cite{liu2017robust}, \cite{kinay2019multi}) to cope with uncertainty. \citet{noyan2012risk} propose risk averse stochastic programming for single-objective problems in disaster management. This work was extended in \cite{noyan2019two} to a multi-criteria optimization approach with a two-stage stochastic programming model. As in our present work, \cite{noyan2019two} use the CVaR to represent risk averseness. However, the multi-criteria decision approach is different from ours, insofar as \cite{noyan2019two} apply the CVaR in the context of multivariate stochastic dominance constraints, whereas we use the concept of Pareto efficiency (the most prominent concept in multi-objective optimization) to determine the trade-off between our considered objective functions.

From a methodological point of view, some of these studies use the $\epsilon$-constraint method \citep{laumanns2006efficient} to find the optimal \textit{Pareto frontier} (e.g., \cite{tricoire2012bi}, \cite{rath2016bi}), while others develop metaheuristic algorithms to approximate the Pareto frontier for large-size instances (e.g., \cite{haghi2017developing}). \textcolor{black}{Let us give details of the mentioned works in the HRC literature, which tackle uncertain multi-objective facility location problems using the concept of Pareto efficiency.} \cite{tricoire2012bi} develop a bi-objective two-stage stochastic model to deal with distribution center selection for relief commodities and delivery planning. The authors consider the demand as a random parameter and approximate them by a sample of randomly generated scenarios. A new solution approach is proposed based on the $\epsilon$-constraint method for the computation of the Pareto frontier of the bi-objective problem using a branch-and-cut algorithm. Finally, they test the proposed algorithm using data of a real-world application in Senegal. In this paper, we address the case of Sengeal using a different model. \cite{rath2016bi} formulate several variants of a two-stage bi-objective stochastic programming model in the response phase of disasters. They apply the $\epsilon$-constraint method to compute the Pareto frontier. They evaluate the benefit of the stochastic bi-objective model in comparison to the deterministic bi-objective model using the value of stochastic solution (VSS) measure. \cite{haghi2017developing} develop a metaheuristic algorithm for a multi-objective location and distribution model with pre/post-disaster budget constraints for goods and casualties logistics. In order to handle the uncertainties, a robust optimization approach is embedded into the $\epsilon$-constraint method. To solve large-size instances, they propose a metaheuristic algorithm that is a combination of a genetic algorithm and simulated annealing \citep{kirkpatrick1983optimization}.  
\cite{hinojosa2014two} investigate a two-stage stochastic transportation problem with uncertain demands, considering an overall objective function that is composed of total cost associated with the selected links in the first decision stage, and expected distribution cost in the second decision stage.\cite{fernandez2019new} deal with a two-stage stochastic mixed-integer model for a fixed-charge transportation problem with uncertain demand on the assumption that the decision maker is risk-averse. Risk aversion is represented by using the CVaR in the objective function. 
\cite{filippi2019single} present a bi-objective facility location problem where the first objective is the minimization of cost. The second objective aims at maximizing fairness by considering a conditional $\beta$-mean measure. The conditional $\beta$-mean is a concept related to the CVaR, but note that in \cite{filippi2019single}, it is applied to fairness quantification rather than to the measurement of risk. The authors develop a weighted-sum method to generate the Pareto-frontier for small/medium-size instances. Besides, a Benders decomposition method is employed to deal with the large-size instances. Comparison of the solution shows the efficiency of their proposed methods. \cite{parragh2020branch} propose an uncertain bi-objective facility location problem considering stochastic demand in a disaster relief context. They formulate the uncertainty using a scenario-based two-stage risk-neutral stochastic approach. The authors integrate the L-shaped method into a bi-objective branch-and-bound framework to deal with the problem. They test and compare different cutting-plane schemes on instances with varying numbers of samples.

\textcolor{black}{\textit{Contributions of the paper.} } The focus of our paper is on the evaluation of different combinations of two criterion space search methods, the well-known $\epsilon$-constraint \citep{laumanns2006efficient} and the recently proposed balanced box \citep{boland2015criterion} methods and three different uncertainty approaches, the widely used risk-neutral two-stage stochastic programming approach, adaptive (two-stage) robust optimization and two-stage risk-averse stochastic programming using CVaR. It aims at finding an efficient method to generate the entire Pareto frontier by considering last-mile setting assumptions. It is assumed that the demand of affected people is uncertain.

In addition, as it is not possible to solve large-size instances to optimality within a reasonable time limit, we propose an iterative mixed integer programming (MIP)-based matheuristic.

\section{Problem description}
\label{sec:problemdec}

Due to scarcity of resources, in many disaster relief situations, including slow-onset disasters or sudden-onset disasters, humanitarian organizations encounter challenging logistical \textit{last-mile} operations \citep{rancourt2015tactical, balcik2008facility}. On one hand, the humanitarian organizations’ goal in such a context is to reduce beneficiaries’ vulnerability by maximizing coverage. On the other hand, they face limited monetary resources and want to reduce their costs. Therefore DMs face two conflicting objectives, and a trade-off solution needs to be identified.

\begin{figure}
	\includegraphics[width=\linewidth]{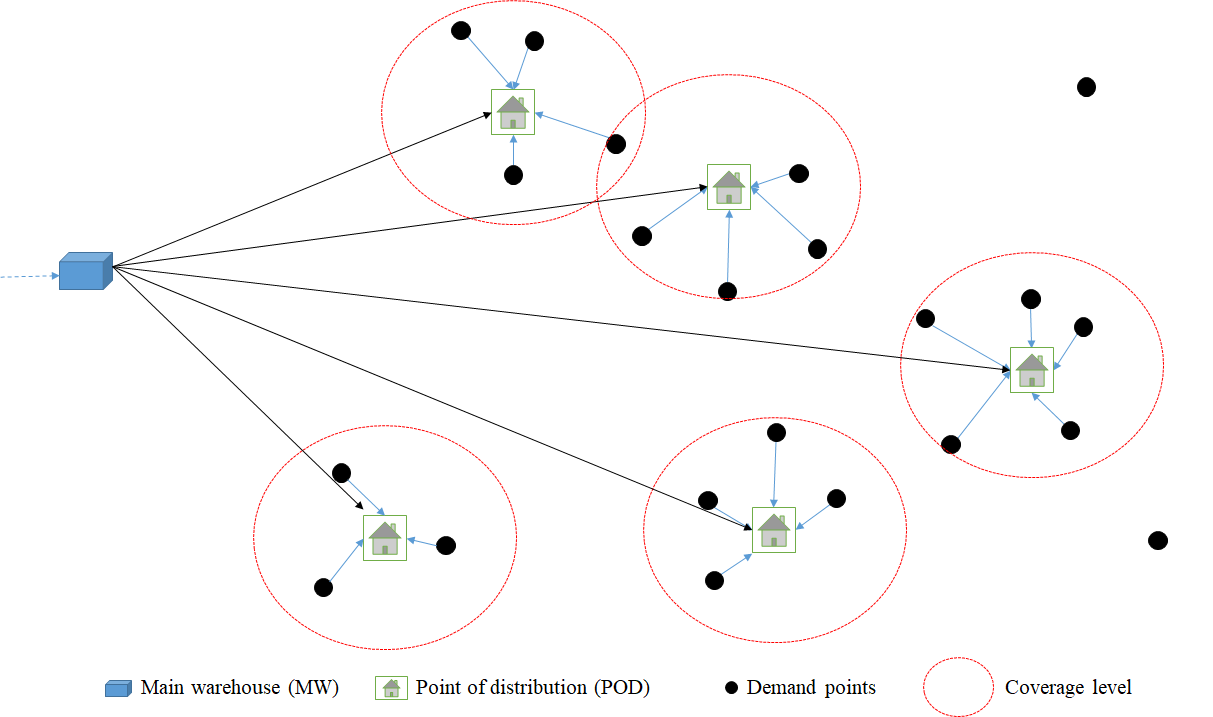}
	\caption{Last Mile Relief Network\\
	{\footnotesize 	Source. Modified from \cite{noyan2016stochastic}}
}
	\label{disas-fig}
\end{figure}

We consider a \textcolor{black}{Bi-SSCFL} \textcolor{black}{with simplified assumptions} of a last-mile network, motivated by the drought case studies presented by \cite{tricoire2012bi} and \cite{rancourt2015tactical}, as well as the earthquake case presented by \cite{noyan2016stochastic} in which the authors design a last-mile aid distribution network. Unlike \cite{tricoire2012bi}, the two latter studies rely on a single objective approach, albeit the authors acknowledge and analyze different objectives. The main framework for the cases is the same (Figure \ref{disas-fig}). \textcolor{black}{Our model is an uncertain bi-objective extension of an uncapacitated deterministic facility location problem proposed by \cite{rancourt2015tactical} to design a real food aid network in Kenya, Sub-Saharan Africa, in which the objective is to minimize the welfare cost of all the involved stakeholders.} The aim is to find the best locations to position PODs in the relief phase of a disaster where beneficiaries walk to these PODs to pick up the aid packages. Two objective functions are considered: the first objective is the minimization of location costs, and the second objective is to minimize the number of uncovered demand of beneficiaries (effectiveness-related criterion). 

In this paper, we consider the amount of demand as an uncertain parameter, as the demand of the affected people depends on the crisis's intensity which is uncertain (single source of uncertainty). Without loss of generality, multiple sources of uncertainty (e.g., the demand of the affected people, the capacity of PODs, the availability of transportation links, etc.) can also be considered in the model.

The problem is formulated on a graph $G= \left(V_0, A\right)$, where $V_0$ is the set of nodes, and $A$ is the set of arcs. $V_0$ can be partitioned into $\{\{0\}, I\}$, where $\{0\}$ is the main warehouse (MW), supplying the selected PODs, and $I$ is the set of affected nodes and a superset of the set of potential PODs, $J$ (without loss of generality we assume the set of potential PODs is a subset of the set of affected nodes  $J\subseteq I$: If a potential POD is not affected, it can be represented as a virtual affected node with zero beneficiaries). We assume that the location of the MW is determined before the optimization.

Another assumption is that beneficiaries in each affected location ($i\in I$) will walk only to the one opened POD ($j\in J$), which they are assigned to, in order to fulfill their demand ($q_i^{(s)}$) (single-source assumption). Moreover, beneficiaries will not walk to any PODs, if their distance ($d_{ij}$) from an opened POD is more than a certain distance threshold ($d_{max}$).

The first stage decision (represented by decision variables $y_j$) is the location of PODs with limited capacity ($c_j$) and given fixed cost ($\gamma_j$). This decision has to be made before the realization of the uncertain parameters. The second stage decisions determine the assignments of demand nodes to the opened PODs (decision variables $x_{ij}^{(s)}$) and the amount of relief items delivered to each POD (decision variables $u_j^{(s)}$), which are determined based on first stage decisions and realized uncertain information. Table \ref{tab1}, summarizes the employed notation.

\begin{table}[H]

		\caption{Notation                                                     }
	\vspace{-10pt}
	\begin{tabular}{llllllll}
		\toprule
		Sets & \multicolumn{7}{c}{} \\
		$I$ & \multicolumn{7}{l}{set of demand nodes} \\
		$J$    & \multicolumn{7}{l}{set of potential PODs, $J\subseteq I$} \\
		$S$   & \multicolumn{7}{l}{set of scenarios} \\
			\multicolumn{8}{c}{} \\
		Certain Parameters & \multicolumn{7}{c}{} \\
		$\gamma_j$ & \multicolumn{7}{l}{fixed cost to open POD $j \in J$} \\
		$c_j$    & \multicolumn{7}{l}{maximum capacity of POD $j \in J$} \\
		$\psi (d_{ij})$   & \multicolumn{7}{l}{coverage level: 1 if $d_{ij} \le d_{max}$, 0 otherwise} \\
			\multicolumn{8}{c}{} \\
		Uncertain Parameters & \multicolumn{7}{c}{} \\
		$q_i^{(s)}$    & \multicolumn{7}{l}{demand of demand point $i \in I$ under scenario $s\in S$} \\
			\multicolumn{8}{c}{} \\
		Decision Variables & \multicolumn{7}{c}{} \\
		$y_j$    & \multicolumn{7}{l}{1 if POD $j \in J$ is selected to be open, 0 otherwise} \\
		$x_{ij}^{(s)}$   & \multicolumn{7}{l}{1 if demand point $i \in I$ is assigned to POD $j \in J$ under scenario $s\in S$, otherwise 0} \\
		$u_j^{(s)}$    & \multicolumn{7}{l}{quantity of supply delivered to POD $j \in J$ under scenario $s\in S$} \\
		\multicolumn{8}{c}{} \\
		Objective Functions & \multicolumn{7}{c}{} \\
		$f_1$    & \multicolumn{7}{l}{first stage objective: Total (location) opening cost} \\
		$f_2$   & \multicolumn{7}{l}{second stage objective: Total amount of uncovered demand (uncertain objective)} \\
		\bottomrule
	\end{tabular}%
	\label{tab1}%
\end{table}%

\subsection{Base formulation}
\label{subsec:base formulation}

In mathematical terms, a nominal formulation of the first stage is the following. Note that the precise meaning of minimizing $f_2$ is not yet specified, we shall return to this issue in Subsection~\ref{subsec:uncertainty}. 

\begin{align}
&\text{min}& &f_1=\sum_{j \in J} \gamma_{j}y_{j}&\\
&\text{min}& &f_2= Q(y,\xi)&\\
\vspace{10pt}&\text{s.t.} & &y_j\in \{0,1\}& &\forall j \in J&
\end{align}

\textcolor{black}{The first objective (1) minimizes opening costs, and the second one (2) minimizes the uncovered demand of the affected people.} Here $\xi$ denotes the random data and $Q(y,\xi)$ is the second objective function associated with the second stage of the problem for a given decision resulting from the first stage. It represents the uncovered demand resulting from decision $y$, if the uncertain data realize as $\xi$. The second stage is formulated as follows for the realization of the random data under scenario $s\in S$ ($\xi^s=(q^s)$):

\begin{align}
&\text{min}& &Q(y,\xi^s)=\sum_{i \in I} q_i^{(s)}-\sum_{j \in J} u_j^{(s)}&\\
&\text{s.t.}& &\sum_{j \in J} \psi(d_{ij})x_{ij}^{(s)}\leq 1& &\forall i \in I&\\
& & &u_j^{(s)}\leq c_j y_j& &\forall j \in J&\\
& & &u_j^{(s)}\leq \sum_{i \in I} q_i^{(s)}\psi(d_{ij})x_{ij}^{(s)}& &\forall j \in J&\\
& & &x_{ij}^{(s)}\leq y_j& &\forall i\in I, \forall j \in J&\\
& & &x_{ij}^{(s)}\in \{0,1\}& &\forall i\in I, \forall j \in J&\\
& & &u_j^{(s)}\in \mathbb{Z^+}& &\forall j \in J&
\end{align}

Constraint (5) indicates coverage constraints. It makes sure that any part of the demand at $i$ is only covered at most once. Constraints (6) ensure that the capacity of POD $j$ is not exceeded. Constraints (7) link the coverage variables with the assignment variables: the covered demand with POD $j$ cannot be higher than the actual demand assigned to this node. Constraints (8) guarantee that a demand node can only be assigned to an opened POD, and finally, constraints (3), (9) and (10) give the domains of the variables.

\subsection{Uncertainty treatment}
\label{subsec:uncertainty}

To deal with parameter uncertainty, as mentioned in the introduction, we use three different approaches, two-stage risk-neutral stochastic programming (expected value), two-stage (adaptive) robust optimization, and two-stage risk-averse stochastic programming (CVaR). We consider a scenario-based framework to make a fair comparison among all the above-mentioned approaches. We characterize the uncertainty in parameters via a finite discrete set of random scenarios ($s\in S$, $S=\{1,...,N\}$) where each scenario has the same probability $\frac{1}{N}$. In the following, we provide the corresponding models.

\subsubsection{Two-stage risk-neutral stochastic optimization}
\label{subsec:twostage}

As mentioned before, the nature of the decision-making process in our problem is sequential: the decision on the location of the POD has to be taken before the realization of uncertain data. The assignment of demand points to an opened POD is done afterwards. We capture this setting through a two-stage stochastic programming model and in a risk-neutral context by taking the expected value of the random variable  $f_2$ as the evaluation of second-stage costs. 
Representing the probability distribution by $N$ equiprobable scenarios as described above, the second-stage objective function ($f_2$) will be replaced by the expected value of the uncovered demand in the stochastic model. So, the gained deterministic counterpart of the bi-objective model (1) to (10) is as follows (\textbf{M1}):

\begin{align}
&\text{min}& &\sum_{j \in J} \gamma_{j}y_{j}&\\
&\text{min}& &E(Q(y,\xi^s))=\frac{1}{N}\sum_{s \in S}(\sum_{i \in I} q_i^{(s)}-\sum_{j \in J} u_j^{(s)})&
\end{align}
$$\hspace{-3cm}\text{s.t.}\hspace{4cm} (5)-(10) \quad \forall s \in S, (3)$$


\subsubsection{Adaptive robust optimization}
\label{subsec:adaptive}

The above two-stage stochastic approach is based on the expected value, whereas minmax robustness is based on the worst-case scenario (pessimistic view). Adaptive robust optimization, in which the recourse function is a worst-case value over a set of scenarios (also called “minimax two-stage stochastic programming”), might give less conservative solutions according to \cite{bertsimas2010models} who compared to static (minmax) robustness. By replacing $f_2$ with its worst value in (1)-(10), the deterministic counterpart of the bi-objective model is as follow (\textbf{M2}):

\begin{align}
&\text{min}& &\sum_{j \in J} \gamma_{j}y_{j}&\\
&\text{min}& &\max_{s\in S}\hspace{1mm}(Q(y,\xi^s)=\sum_{i \in I} q_i^{(s)}-\sum_{j \in J} u_j^{(s)})&
\end{align}
$$\hspace{-3cm}\text{s.t.}\hspace{4cm} (5)-(10) \quad \forall s \in S, (3)$$

\subsubsection{Two-stage risk-averse stochastic optimization}
\label{subsec:robust}

So far, to cope with randomness in our problem, we have employed the expected value of coverage, which corresponds to a risk-neutral approach, and the worst-case value of the coverage, which is the adaptive robustness approach. These two approaches are at the two ends of the risk averseness spectrum, in which the risk-neutral approach gives a comparably optimistic view  whereas the adaptive robust approach represents a pessimistic view of the random outcome. Sometimes, DMs are not risk-neutral, but less risk-averse than assumed by the adaptive robust approach. For these cases, we consider a two-stage risk-averse stochastic optimization model where the degree of risk-aversion can be specified. We use the widely applied conditional value-at-risk (CVaR) as the risk measure in our study which leads to a computationally tractable model. 

In this paper we employ classical linear representation of CVaR that is valid on the assumption of an underlying finite discrete probability space. \textcolor{black}{The general concept of CVaR and its classical representation is briefly discussed in Appendix~\ref{sec:AppA}}.

\paragraph{The classical representation of CVaR  (CCVaR)}

In our problem, we calculate the CVaR value for the second-stage objective $Q(y,\xi^s)$ with $s=1,...,N$. 
The deterministic equivalent formulation of our model is as follows (\textbf{M3}): 

\begin{align}
&\text{min}& &\sum_{j \in J} \gamma_{j}y_{j}&\\
&\text{min}& &CVaR_\alpha(Q(y,\xi^s))=\eta+\frac{1}{1-\alpha}\frac{1}{N}\sum_{s \in S}  w_s&\\
&\text{s.t.}& &w_s\geq Q(y,\xi^s)-\eta& &\forall s \in S&\\
& & &w_s\geq 0& &\forall s \in S&\\
& & &\eta \in \mathbb{R}& &&
\end{align}
$$\hspace{-4cm} (5)-(10) \quad \forall s \in S, (3)$$

This formulation is equivalent to the two-stage stochastic (expected value) approach in the special case $\alpha = 0$, and it is equivalent to the adaptive (two-stage) robust approach for sufficiently large value of $\alpha,$ $\alpha \to 1$.

\section{Solution approach}
\label{sec:method}

In this section, we first describe two exact bi-objective frameworks that are used to solve small instances. Thereafter, we propose a matheuristic method to deal with large-size instances. It uses the MIP as a backbone and a local search algorithm to heuristically generate better solutions. \textcolor{black}{For the definition of solution concepts in a multi-objective setting, we refer to Appendix~\ref{sec:AppB}}.

\subsection{Multi-objective exact solution techniques}
\label{subsec:multiobj}

Most of the papers in the literature of HRC have deployed metaheuristic techniques to solve large-size (i.e., real-world) instances of mathematical models. Metaheuristic techniques cannot provide performance guarantees, whereas exact methods can. Among exact methods, criterion space search methods (i.e., methods that search in the space of objective function values) are often computationally more efficient than decision space search methods (i.e., methods that search in the space of feasible solutions) \citep{boland2015criterion}.  To find a good approximation of the entire Pareto frontier, we employ two criterion space search frameworks, the well-known $\epsilon$-constraint method due to its simplicity \citep{laumanns2006efficient} and the recently developed balanced box method \citep{boland2015criterion} since it has been shown that it performed well for large-size problems. It is worth noting that these frameworks are exact algorithms that give an approximation of the Pareto frontier if they are terminated before completion, e.g. if a tight time limit is employed.  

The $\epsilon$-constraint method generates all NDPs step by step iteratively. Let us describe the method for the special case \textcolor{black}{of two objectives}. It starts by finding the endpoints of the Pareto frontier using lexicographic optimization. Lexicographic optimization is performed as follows: we optimize the first objective function ($f_1$) and gain the optimal minimum value $f_1=f^{min}_1$. Then, we optimize the second objective function ($f_2$) by adding constraint $f_1=f^{min}_1$ to the model in order to keep the optimal solution value of the first optimization and $f^{max}_2$ is obtained. The solution to this procedure gives the first extreme point ($Z^T$). This procedure can be represented in the following notation (note that $R(z^1,z^2)$ is the box in the criterion space defined by the points $z^1=(f_{1}^1,f_{2}^1)$ and $z^2=(f_{1}^2,f_{2}^2)$ where $f_{1}^1 \leq f_{1}^2$ and $f_{2}^2 \leq f_{2}^1$):

\begin{align}
Z^T=lexmin\{f_1,f_2: f \in R((-\infty,\infty),(\infty,-\infty))\}
\end{align}

The second extreme point ($Z^B$) is found by optimizing the second objective function ($f_2$) and obtaining the optimal minimum value $f_2=f^{min}_2$. Then, the first objective function ($f_1$) is optimized by adding constraint $f_2=f^{min}_2$ to the model and $f^{max}_1$ is gained:

\begin{align}
Z^B=lexmin\{f_2,f_1: f \in R((-\infty,\infty),(\infty,-\infty))\}
\end{align}

The $\epsilon$-constraint algorithm enumerates each NDP one by one from the direction of one endpoint to the other. In all iterations, one and the same of the objective functions is considered as the main objective and the other in terms of an $\epsilon$-constraint. The parameter $\epsilon$ is set as given in lines 5 and 7 in Algorithm \ref{alg1}.

We note that the choice of direction in the $\epsilon$-constraint method can play an essential role in the performance of the algorithm. In our case, we optimize $f_2$ as the main objective and handle $f_1$ in the $\epsilon$-constraint based on preliminary results. That is, the Pareto solutions are enumerated from $Z^B$ to $Z^T$. Since the first objective takes integer values in our problem, the smallest difference value ($d$) is $1$. Algorithm \ref{alg1} shows the procedure of the $\epsilon$-constraint framework.

\begin{algorithm}
	\caption{$\epsilon$-constraint framework}
	\begin{algorithmic}[1]
		\STATE \textbf{input} : $Z^T=(f^{min}_1,f^{max}_2), Z^B=(f^{max}_1,f^{min}_2)$
		\STATE $L \gets \emptyset$  
		\STATE $\epsilon \gets f^{max}_1-d$
		\REPEAT
		\STATE $x \gets lexmin\{f_2,f_1: f_1 \leq \epsilon\}$
		\STATE $L \gets L \cup x$
		\STATE $\epsilon \gets f_1(x)-d$
		\UNTIL $\epsilon \geq f^{min}_1$
		\RETURN $L$
	\end{algorithmic}
	\label{alg1}
\end{algorithm}

On the other hand, the balanced box method keeps a priority queue of boxes in criterion space in non-decreasing order of their areas. In the beginning, the priority queue is empty. So, similar to the $\epsilon$-constraint method, the algorithm first finds the endpoints of the Pareto frontier ($Z^T, Z^B$). These two points define the initial box $R(Z^T, Z^B)$. It contains all not yet found NDPs. Subsequently, all the other boxes are generated and explored iteratively.

In each iteration, the largest box in the priority queue pops out, and the algorithm then splits the box horizontally into two equal parts $R^T$ and $R^B$. It first explores the bottom box for a NDP by optimizing the first objective function. It then explores the top box for a NDP with minimization of the second objective function. The procedure can be repeated until no more boxes are left unexplored in the priority queue. The pseudocode of this framework is presented in Algorithm \ref{alg2}.

\begin{algorithm}
	\caption{Balanced box framework}
	\begin{algorithmic}[1]
	\STATE \textbf{input} : $Z^T=(f^{min}_1,f^{max}_2), Z^B=(f^{max}_1,f^{min}_2)$
	\STATE $L \gets \emptyset$ 
	\STATE $P \gets R(Z^T,Z^B)$   
	\WHILE {$P \neq \emptyset$}
	\STATE $R^B \gets R((f^{min}_1,(\frac{f^{min}_2+f^{max}_2}{2})),Z^B)$
	\STATE $\bar{Z}^T \gets lexmin(f_1,f_2,f \in R^B)$

	\IF { $\bar{Z}^T \neq Z^B$}
\STATE $L \gets L \cup \bar{Z}^T$ 
\STATE $P \gets P \cup R(\bar{Z}^T,Z^B)$
\ENDIF
\STATE $R^T \gets R(Z^T,(\bar{Z}^T-d)$
\STATE $\bar{Z}^B \gets lexmin(f_2,f_1,f \in R^T)$
\IF { $\bar{Z}^B \neq Z^T$}
\STATE $L \gets L \cup \bar{Z}^B$ 
\STATE $P \gets P\cup R(Z^T,\bar{Z}^B)$
\ENDIF
	\ENDWHILE
	\RETURN $L$
	\end{algorithmic}
	\label{alg2}
\end{algorithm}

\subsection{A matheuristic solution algorithm}

Although small instances can be solved to optimality by the exact solution techniques, we develop a matheuristic method to compute high-quality approximated Pareto frontiers for the large-size instances. It combines the feature of the $\epsilon$-constraint method and \textcolor{black}{bi-objective generalizations of two single-objective variable fixing heuristic algorithms, namely local branching \citep{fischetti2003local} and relaxation induced neighborhood search (RINS) \citep{danna2005exploring}}

Local branching \citep{fischetti2003local} is a general MIP-based framework which explores a neighborhood of a feasible reference solution $\bar{y}$ to identify a better solution. Assume $\bar{A}=\{j \in V: \bar{y}_j=1\}$ is the set of indices of potential POD variables where their value is equal to one. To construct an $l$-OPT neighborhood $N(\bar{y},l)$ of the reference solution, the following local branching constraint is added to the model:

\begin{align}
\sum_{j \in \bar{A}}(1-y_j)+\sum_{j \in J \setminus \bar{A}}y_j \leq l
\end{align}

Constraint (22) states that at most $l$ variables of the reference solution $\bar{y}$ switch their values either from 1 to 0 or from 0 to 1.

\textcolor{black}{Note} that the cardinality of the binary support of any feasible solution is a constant in our model. Therefore, constraint (22) can equivalently be written in its \textit{asymmetric} form:

\begin{align}
\sum_{j \in \bar{A}}(1-y_j) \leq l^\prime
\end{align}

\textcolor{black}{Moreover, RINS \citep{danna2005exploring} is another general MIP-based mechanism to search a neighborhood of a feasible solution. It uses the information retrieved within the branch-and-bound tree comparing the incumbent solution and the continuous relaxation solution concerning variables which take the same values in both solutions. By fixing the variables with identical values, its focus is on the variables with different values. Similarly, comparing the Pareto solutions corresponding to NDPs along the Pareto frontier in the Bi-SSCFL problem, with identical opening costs, show that they have many variables with the same values. We observed that when moving along the frontier solution by solution from one end-point ($Z^B$) to the other one ($Z^T$), very often one of the opened facilities (i.e., $y_j=1$) is closed, and the set of closed facilities is kept unchanged (i.e., $y_j=0$). Therefore the variables with value $0$ could be fixed.}

\textcolor{black}{Assume $A=\{j \in J: y_{j}^a=y_{j}^b=0\}$ is the set of potential POD variables related to neighbor solutions $y^a$ and $y^b$, where their value is equal to zero.}
\textcolor{black}{We developed the following procedure in order to find a high-quality approximation of the Pareto frontier (cf. Algorithm~\ref{alg3}): \\
(i) we find the optimal Pareto efficient endpoints ($Z^T$ and $Z^B$) and create the initial box $R(Z^T,Z^B)$; 
(ii) noting that we move along the frontier in the direction from $Z^B$ to $Z^T$ in our $\epsilon$-constraint method, we add the linear constraint of $y_{j}^b=0\quad \forall j\in A$, based on solution $Z^B$, and we approximate the next efficient solution, namely $y^a$, if the MIP problem is solved in a certain time limit (TiLim). Otherwise, we approximate the efficient solution $y^a$ by solving the linear relaxation (LR) of the model, rounding the fractional values in order to obtain an integer solution, and we deploy local branching to approximate a better solution; (iii) we generate the entire Pareto frontier iteratively updating the $y^b \gets y^a$.}
	
\textcolor{black}{\cite{leitner2016ilp} propose a two-phase heuristic approach that uses a generalization of RINS called BINS and local branching to solve bi-objective binary integer linear programs. Our method differs from theirs in two ways. First, we combine an adaptation of RINS with the $\epsilon$-constraint framework to explore the criterion space, whereas \cite{leitner2016ilp} apply BINS in the first phase of their method, to explore each box found in a weighted-sum framework. Secondly, our method employs the local branching scheme where the problem cannot be solved to optimality fast enough by only fixing the variables, they applied local branching to refine the NDPs found in the first phase.}
	
	\floatname{algorithm}{\color{black}Algorithm}
	\begin{algorithm}
		\caption{\textcolor{black}{Matheuristic framework}}
		\color{black}
		\begin{algorithmic}[1]
			\STATE \textbf{input} : $Z^T=(f^{min}_1,f^{max}_2), Z^B=(f^{max}_1,f^{min}_2)$
			\STATE $L \gets \emptyset$  
			\STATE $\epsilon \gets f^{max}_1-d$
			\REPEAT
			\IF{TiLim is not reached}
			\STATE $x \gets lexmin\{f_2,f_1: f_1 \leq \epsilon, y_{j}^b=0\quad \forall j\in A\}$
			\ELSE
			\STATE $x \gets LR\{f_2: f_1 \leq \epsilon, y_{j}^b=0\quad \forall j\in A\}$
			\STATE $LocalBranching(f_1(x))$
			\ENDIF
			\STATE $L \gets L \cup x$
			\STATE $\epsilon \gets f_1(x)-d$
			\STATE $y^b \gets y^a$
			\UNTIL $\epsilon \geq f^{min}_1$
			\RETURN $L$
		\end{algorithmic}
		\label{alg3}
\end{algorithm}

\textcolor{black}{It is worth mentioning that based on our preliminary computational experiments, computing the endpoints $Z^T$ and $Z^B$ was not expensive in our problem. However, in cases where the exact computation of the two endpoints is too demanding, we suggest to use the local branching approach to approximate them.}

\section{Computational study}
\label{result}
In the following, we present the computational experiments conducted to evaluate the performance of the proposed methods. We first describe the test instances and then, the quality indicators employed to assess the performance of the methods. This is followed by a recapitulation of the expected value of perfect information (EVPI), and the value of stochastic solution (VSS) used to show the value of incorporating risk measurement. Thereafter, we discuss the computational results.

All experiments have been implemented in C++ using the Concert Technology component library of IBM ®ILOG ®CPLEX ®12.9 as a MIP solver where multi-threading is disabled. \textcolor{black}{The algorithms are run on a cluster where each node consists of two Intel Xeon X5570 CPUs at 2.93 GHz and 8 cores, with 48GB RAM.} A time limit (TL) of 7200 s has been applied. \textcolor{black}{In addition, the TiLim parameter for the matheuristic technique is set to 150 s}.

\subsection{Test instances}
\label{test instances}
To test and evaluate the integration of the above-mentioned uncertainty approaches into the stated multi-objective frameworks, a data set inspired by a real-world study for a slow-onset drought case presented in \citet{tricoire2012bi} \textcolor{black}{has been utilized}. \citet{tricoire2012bi} provided \textcolor{black}{us with the} data set from the region of Thies in western Senegal.  Senegal is a developing country in sub-Saharan Africa where droughts occur frequently. Politically, the region is split into 32 rural areas where each contains between 9 and 31 villages (demand nodes). In total, the region contains 500 nodes. In order to test the performance of the proposed solution approaches, we create new networks and artificially generate larger instances with between 21 and 500 nodes by aggregating the demand nodes in the region. For example, the first instance contains 21 nodes, and the second one has 44 nodes, i.e., an instance with 23 nodes is added to the first 21 nodes, \textcolor{black}{and the second instance of 44 nodes is derived}. Our data set contains 23 test instances.

The distance matrix is obtained by road distances between each pair of nodes. Opening costs for PODs are assumed to be identical for all demand locations (5000 cost units). We also assume that the affected people in a demand node will walk to the closest POD if the distance is less than 6km. \textcolor{black}{The capacity of each POD is set equal to the population size of the demand nodes times three.} We start by considering 10 sample scenarios ($|S|$= 10) to cope with demand uncertainty. Thereafter, to assess the algorithms for a larger number of scenarios, samples with 100, 500, and 1000 scenarios are generated. \textcolor{black}{The uncertain demand of each demand node (village) $i$ is considered as $\xi_iq_i$, where $\xi_i$ and $q_i$ are the \textit{uncertainty factor} and the population size of the node $i$, respectively. The uncertainty factor is a sum of a random baseline term ($\xi_{base}$), which is the same for all nodes, and a correction term, which is unique for each demand point. Correction terms for different demand points are assumed as stochastically independent. In mathematical formulas, $\xi_i= \xi_{base} - \lambda_2 + 2\lambda_2Z_i$, where $\xi_{base}=\bar{\xi}-\lambda_1+2\lambda_1Z$. Therein, $\bar{\xi}$ is a constant value set to 1, $Z$ and $Z_i$ are random numbers uniformly distributed between 0 and 1, and $\lambda_1=\lambda_2\geq 0$ are constant parameters chosen to be 0.5.} We refer to \citet{tricoire2012bi} for further details.

\subsection{Solution quality indicators}
Different quality indicators for multi-objective approximation algorithms have been proposed in the literature. \citet{zitzler2003performance} review the existing quality measures. We employ the hypervolume indicator to assess the performance of the proposed methods.

Given a bi-objective framework with minimization objective functions, let $X$ be the set of feasible solutions in decision space and $Z=f(X)$ be the set of feasible solutions in criterion space. We assume $A \subseteq Z$ is an approximation set of a Pareto frontier, and $R \subseteq Z$ is a reference set, e.g., the true Pareto frontier.

The hypervolume indicator ($I_H$) \citep{zitzler1998multiobjective} computes the area of the portion of the criterion space, which is weakly dominated by approximation set $A$ with respect to a reference point (e.g., the Nadir point). The complete Pareto frontier generally has a maximum value of $I_H$. The higher the value of $I_H$, the higher is the quality of the approximated Pareto frontier. 

\subsection{Value of using a risk measure stochastic model}
\label{EVPI and VSS}

Due to the computational challenges of decision making under uncertainty, one of the first questions that DMs might come up with is whether it pays off to consider a stochastic instead of a deterministic model. An answer can be given by using two well-known measures: the expected value of perfect information (EVPI) and the value of stochastic solution (VSS). These measures assess the value of using a single objective two-stage risk-neutral programming model (see, \citet{birge2011introduction}).

In this section, we start by recalling these two measures. After that, similarly as in \citet{noyan2012risk} where the measures are extended to a single-objective two-stage mean-risk model involving the CVaR, we adapt these concepts to our bi-objective two-stage CVaR model.

The EVPI measures the expected gain of perfect information over the stochastic solution. It is the difference between the wait-and-see solution (WS) and the solution obtained by solving the risk-neutral stochastic model referred to as the recourse problem (RP). WS is obtained by solving the model for each scenario as realized data with a probability of 1, and then taking the expected objective function value:

\begin{align}
\text{WS}= E (Q(\bar{y}(\xi^s),\xi^s))
\end{align} 

Therein, $\bar{y}(\xi^s)$ denotes the optimal solution of the individual problem for each scenario. Then, by definition: EVPI = RP$-$WS. 

On the other hand, the VSS evaluates the stochastic model in comparison to an expected value solution, where the latter is defined as the solution of the deterministic problem obtained by replacing each random parameter by its expected value. The VSS is calculated as the difference between the RP solution value and the solution value EEV of the expected value problem. EEV is obtained by solving a deterministic expected value scenario problem in the first step. Then, the resulting optimal first-stage variables are saved and fixed in the model, and the second stage is solved:

\begin{align}
\text{EEV}= E(Q(\bar{y}(\bar{\xi}),\xi^s))
\end{align} 

Where, $\bar{\xi}= E (\xi^s)$, and $\bar{y}(\bar{\xi})$ is the the expected value solution. Then, VSS is calculated as follow: VSS = EEV$-$RP. 

As it is mentioned, these two measures are based on expected values which are used to assess a risk-neutral stochastic model. However, they cannot be used directly to evaluate a two-stage risk-averse stochastic model. Therefore, we adapt these measures and apply the risk function (CVaR) instead of the expected value in computing the WS and EEV approach. More precisely, we use the following measures for the two-stage CVaR model at the risk level of $\alpha$ (see, \cite{noyan2012risk}):

\begin{align}
\text{RVPI}(\alpha)= \text{RRP}(\alpha) - \text{RWS}(\alpha)\\
\text{RVSS}(\alpha)= \text{REV}(\alpha) - \text{RRP}(\alpha)
\end{align} 

Therein, RRP$(\alpha)$ is the solution obtained by solving the risk-averse CVaR model  at the risk level of $\alpha$, RWS$(\alpha)= \text{CVaR}_\alpha (Q(\bar{y}(\xi^s),\xi^s))$ is obtained by solving the model for each scenario as realized data with a probability of 1 and then taking the CVaR objective function value at the risk level of $\alpha$, and REV$(\alpha)= \text{CVaR}_\alpha (Q(\bar{y}(\bar{\xi}),\xi^s))$ is obtained by solving a deterministic expected value scenario problem in the first step. Then, the resulting optimal first-stage variables are saved and fixed in the model, and the second stage of the risk-averse model is solved  at the risk level of $\alpha$.

RVPI measures the gain of perfect information based on the CVaR value of the objective values obtained from WS solutions. The RVSS measures the gain from solving the risk-averse model with a specific risk preference. The higher the values of RVSS, the more is the value-added of considering a risk-averse model instead of a risk-neutral problem.

In order to employ these measures to our bi-objective model, we apply them to each non-dominated solution where $f_1$ is bounded as $\epsilon$-constraint, and $f_2$ is the main objective.

\subsection{Results}
In this section, we conduct numerical experiments to assess the proposed approaches. We compare all combinations; $\epsilon$-constraint (e), balanced box (BB) frameworks, and proposed matheuristic (Mat) integrated into deterministic counterparts of the M1, M2, and M3 models. \textcolor{black}{Additionally, we consider the generic feasibility pump based heuristic (FPBH) method recently proposed by \citet{pal2019fpbh}. To implement it, we employed the Julia package, namely "FPBHCPLEX.jl", which uses CPLEX 12.7 as a MIP solver. It is available as an open-source package on GitHub. Since, it was difficult to setup this package on the cluster, FPBH is run on a local computer with Intel® Core™ i5-7200U CPU at 2.50GHz with 16GB RAM. For the sake of a fair comparison, our proposed matheuristic is also run on the same system using CPLEX 12.7 as the solver.} Let \{e-M1, BB-M1, e-M2, BB-M2, e-M3, BB-M3, Mat-M3, \textcolor{black}{FPBH-M3}\} denote the set of all methods considered in this study. Since special cases of M3 model are identical with M1 and M2 models, we just address the performance of Mat method integrated with the M3 model.

We first study the combinations of exact multi-objective methods with two widely used uncertainty approaches, two-stage stochastic, and two-stage robust approach (e-M1, BB-M1, e-M2, BB-M2). Then, we compare the performance of their deterministic counterpart with two-stage CVaR deterministic counterpart in its two special cases where $\alpha= 0$ (risk-neutral) and $\alpha$= $0.9$ (where $|S|$=10 $\to$ the worst value of $\alpha=1-\frac{1}{10}$=0.9). Later, we address the combination of exact and matheuristic multi-objective methods to two-stage CVaR for different levels of risk (e-M3, BB-M3, Mat-M3).

Table~\ref{tab2} indicates the run time of {e-M1, BB-M1, e-M2, BB-M2, e-M3, BB-M3} in seconds for different instances. “TL” indicates that the run-time limit, which we fixed at 7200 s, is reached.   

As can be seen from Table~\ref{tab2} the combination of the $\epsilon$-constraint framework with M1 model (e-M1) finds the complete Pareto frontier even for rather large instances. Furthermore, the deterministic counterpart of two-stage risk-neutral stochastic (M1) and adaptive robust (M2) models outperform the deterministic counterpart of the CVaR model in its special cases ($\alpha$=0 and $\alpha$=0.9). 

\begin{table}[H]
	\centering
	\caption{Run time comparison for M1, M2 and their equivalent special cases of M3 model with $\alpha$=0 and $\alpha$=0.9 combined with exact algorithms for test instances of size 21-500 with sample size 10}
	\small\addtolength{\tabcolsep}{-3pt}
	\makebox[\linewidth]{
		\scalebox{0.85}{
	\begin{tabular}{c|ccccrcccc}
		\multicolumn{1}{r}{} & \multicolumn{9}{l}{$\alpha$ value} \\
		\cmidrule{2-10}    \multicolumn{1}{c}{} & \multicolumn{4}{c}{\textbf{0}} &       & \multicolumn{4}{c}{\textbf{0.9}} \\
		\cmidrule{1-5}\cmidrule{7-10}    \multicolumn{1}{c}{\textit{\#Node}} & \textbf{e-M1} & \textbf{BB-M1} & \textbf{e-M3} & \textbf{BB-M3} &       & \textbf{e-M2} & \textbf{BB-M2} & \textbf{e-M3} & \textbf{BB-M3} \\
		\cmidrule{1-5}\cmidrule{7-10}    \textbf{21} & 1     & 1     & 1     & 2     &       & 2     & 2     & 2     & 2 \\
		\textbf{44} & 16    & 16    & 16    & 22    &       & 21    & 22    & 20    & 26 \\
		\textbf{56} & 18    & 20    & 24    & 31    &       & 20    & 23    & 19    & 33 \\
		\textbf{72} & 28    & 29    & 36    & 51    &       & 61    & 68    & 64    & 103 \\
		\textbf{90} & 56    & 59    & 62    & 110   &       & 98    & 110   & 99    & 179 \\
		\textbf{106} & 115   & 120   & 123   & 229   &       & 226   & 290   & 210   & 380 \\
		\textbf{120} & 128   & 132   & 143   & 241   &       & 139   & 175   & 138   & 262 \\
		\textbf{163} & 421   & 442   & 445   & 818   &       & 769   & 910   & 636   & 1365 \\
		\textbf{182} & 522   & 559   & 577   &  TL     &       & 1555  & 1970  & 1326  & 2548 \\
		\textbf{203} & 820   & 1050  & 956   & 2078  &       & 541   & 652   & 515   & TL \\
		\textbf{254} & 5821  & TL    & 6103  & TL    &       & TL    & TL    & TL    & TL \\
		\textbf{264} & 1950  & 2590  & 3903  & 6800  &       & TL    & TL    & TL    & TL \\
		\textbf{275} & 5499  & TL    & 7133  & TL    &       & TL    & TL    & TL    & TL \\
		\textbf{295} & TL    & TL    & TL    & TL    &       & TL    & TL    & TL    & TL \\
		\textbf{326} & TL    & TL    & TL    & TL    &       & TL    & TL    & TL    & TL \\
		\textbf{355} & TL    & TL    & TL    & TL    &       & TL    & TL    & TL    & TL \\
		\textbf{388} & TL    & TL    & TL    & TL    &       & TL    & TL    & TL    & TL \\
		\textbf{410} & TL    & TL    & TL    & TL    &       & TL    & TL    & TL    & TL \\
		\textbf{436} & TL    & TL    & TL    & TL    &       & TL    & TL    & TL    & TL \\
		\textbf{449} & TL    & TL    & TL    & TL    &       & TL    & TL    & TL    & TL \\
		\textbf{472} & TL    & TL    & TL    & TL    &       & TL    & TL    & TL    & TL \\
		\textbf{482} & TL    & TL    & TL    & TL    &       & TL    & TL    & TL    & TL \\
		\textbf{500} & TL    & TL    & TL    & TL    &       & TL    & TL    & TL    & TL \\
		\bottomrule
		\bottomrule
	\end{tabular}%
}
}
	\label{tab2}%
\end{table}%

Next, we show the computational results for the M3 model \textcolor{black}{on the instances solved in Table~\ref{tab2},} obtained on different levels of $\alpha$ integrated into two exact multi-objective frameworks (Table~\ref{tab3}). These results also indicate that a combination of different $\alpha$ levels with the $\epsilon$-constraint method performs better than their combination with the BB technique.

\begin{table}[htbp]
\centering
\caption{Run time comparison for M3 model combined with exact algorithms and different level of $\alpha$ for test instances of size 21-500 with sample size 10}
\small\addtolength{\tabcolsep}{-3pt}
\makebox[\linewidth]{
	\scalebox{0.85}{
	\begin{tabular}{c|cccccccccccccc}
		\multicolumn{1}{r}{} & \multicolumn{14}{l}{$\alpha$ value} \\
		\cmidrule{2-15}    \multicolumn{1}{c}{} & \multicolumn{2}{c}{\textbf{0}} &       & \multicolumn{2}{c}{\textbf{0.2}} &       & \multicolumn{2}{c}{\textbf{0.5}} &       & \multicolumn{2}{c}{\textbf{0.7}} &       & \multicolumn{2}{c}{\textbf{0.9}} \\
		\cmidrule{1-3}\cmidrule{5-6}\cmidrule{8-9}\cmidrule{11-12}\cmidrule{14-15}    \textit{\#Node} & \textbf{e-M3} & \textbf{BB-M3} &       & \textbf{e-M3} & \textbf{BB-M3} &       & \textbf{e-M3} & \textbf{BB-M3} &       & \textbf{e-M3} & \textbf{BB-M3} &       & \textbf{e-M3} & \textbf{BB-M3} \\
		\cmidrule{1-3}\cmidrule{5-6}\cmidrule{8-9}\cmidrule{11-12}\cmidrule{14-15}    \textbf{21} & 1     & 2     &       & 1     & 2     &       & 2     & 2     &       & 2     & 1     &       & 2     & 2 \\
		\textbf{44} & 16    & 22    &       & 17    & 23    &       & 18    & 26    &       & 17    & 23    &       & 20    & 26 \\
		\textbf{56} & 24    & 31    &       & 22    & 40    &       & 21    & 40    &       & 18    & 22    &       & 19    & 33 \\
		\textbf{72} & 36    & 51    &       & 37    & 54    &       & 42    & 59    &       & 45    & 48    &       & 64    & 103 \\
		\textbf{90} & 62    & 110   &       & 66    & 119   &       & 70    & 135   &       & 73    & 78    &       & 99    & 179 \\
		\textbf{106} & 123   & 229   &       & 120   & 238   &       & 128   & 239   &       & 151   & 161   &       & 210   & 380 \\
		\textbf{120} & 143   & 241   &       & 136   & 263   &       & 143   & 260   &       & 144   & 158   &       & 138   & 262 \\
		\textbf{163} & 445   & 818   &       & 428   & 1028  &       & 461   & 936   &       & 417   & 506   &       & 636   & 1365 \\
		\textbf{182} & 577   & 1176  &       & 579   & 1362  &       & 770   & 1673  &       & 932   & 1053  &       & 1326  & 2548 \\
		\textbf{203} & 956   & 2078  &       & 762   & 3881  &       & 903   & TL    &       & 616   & 699   &       & 515   & TL \\
		\textbf{254} & 6103  & TL    &       & TL    & TL    &       & TL    & TL    &       & TL    & TL    &       & TL    & TL \\
		\textbf{264} & 3903  & 6800  &       & 3422  & TL    &       & 4268  & TL    &       & 6453  & TL    &       & TL    & TL \\
		\textbf{275} & 7133  & TL    &       & 6214  & TL    &       & 6535  & TL    &       & TL    & TL    &       & TL    & TL \\
		\bottomrule
		\bottomrule
	\end{tabular}%
}
}
	\label{tab3}%
\end{table}%

After analyzing the two exact multi-objective techniques, we apply the proposed matheuristic method to the M3 model in order to solve the test instances. As indicated in Table~\ref{tab5}, contrary to the exact methods, the Mat method successfully solves all the instances within the given time-limit and finds an approximation of the entire Pareto frontier. Additionally, we solve a few instances with a larger number of scenarios. The results show that the Mat method also works when the number of scenarios increases (\textcolor{black}{see, Appendix~\ref{sec:AppC}}).

Computational results of the performance measurements are summarized in Table~\ref{tab7}. Table~\ref{tab7} details the performance of the $\epsilon$-constraint method, the BB method, and the Mat method on all of the instances. We report the number of found NDPs and the $I_H$ values. For this purpose, we compute the nadir point as a reference point by calculating the worst objective values over the found optimal efficient set in the $\epsilon$-constraint method. As can be seen from Table~\ref{tab7}, the $\epsilon$-constraint and the BB method do not perform well on larger instances within the time limit. However, the BB method, which bidirectionally explores the criterion space, finds more NDPs. Furthermore, the number of found NDPs and $I_H$ values associated with the Mat method shows that it outperforms the two exact methods on the large-size instances.

\begin{table}[!h]
	\centering
	\caption{Run time comparison of a combination of $\epsilon$-constraint method (e), and the proposed matheuristic (Mat) method with M3 for different $\alpha$/$k$ values for test instances of size 21-500 with sample size of 10 scenarios.}
		\small\addtolength{\tabcolsep}{-3pt}
	\makebox[\linewidth]{
		\scalebox{0.95}{
	\begin{tabular}{c|cccccccc}
		\multicolumn{1}{r}{} & \multicolumn{5}{l}{$\alpha$/$k$ value} &       &       &  \\
		\cmidrule{2-9}    \multicolumn{1}{c}{} & \multicolumn{2}{c}{\textbf{0/10}} &       & \multicolumn{2}{c}{\textbf{0.7/3}} &       & \multicolumn{2}{c}{\textbf{0.9/1}} \\
		\cmidrule{1-3}\cmidrule{5-6}\cmidrule{8-9}    \multicolumn{1}{c}{\textit{\#Node}} & \textbf{e-M3} & \textbf{Mat-M3} &       & \textbf{e-M3} & \textbf{Mat-M3} &       & \textbf{e-M3} & \textbf{Mat-M3} \\
		\midrule
		\midrule
		\textbf{21} & 1     & 0.51  &       & 2     & 0.43  &       & 2     & 0.58 \\
		\textbf{44} & 16    & 3     &       & 17    & 3     &       & 20    & 3 \\
		\textbf{56} & 24    & 4     &       & 18    & 4     &       & 19    & 5 \\
		\textbf{72} & 36    & 11    &       & 45    & 9     &       & 64    & 11 \\
		\textbf{90} & 62    & 18    &       & 73    & 23    &       & 99    & 21 \\
		\textbf{106} & 123   & 32    &       & 151   & 36    &       & 210   & 44 \\
		\textbf{120} & 143   & 37    &       & 144   & 39    &       & 138   & 43 \\
		\textbf{163} & 445   & 183   &       & 417   & 132   &       & 636   & 167 \\
		\textbf{182} & 577   & 212   &       & 932   & 178   &       & 1326  & 172 \\
		\textbf{203} & 956   & 206   &       & 616   & 217   &       & 515   & 210 \\
		\textbf{254} & 6103  & 414   &       & TL    & 589   &       & TL    & 568 \\
		\textbf{264} & 3903  & 479   &       & 6453  & 470   &       & TL    & 478 \\
		\textbf{275} & 7133  & 686   &       & TL    & 781   &       & TL    & 819 \\
		\textbf{295} & TL    & 726   &       & TL    & 729   &       & TL    & 836 \\
		\textbf{326} & TL    & 893   &       & TL    & 2407  &       & TL    & 1182 \\
		\textbf{355} & TL    & 1447  &       & TL    & 1499  &       & TL    & 1619 \\
		\textbf{388} & TL    & 1611  &       & TL    & 1835  &       & TL    & 1931 \\
		\textbf{410} & TL    & 2133  &       & TL    & 4062  &       & TL    & 2847 \\
		\textbf{436} & TL    & 2439  &       & TL    & 2531  &       & TL    & 3608 \\
		\textbf{449} & TL    & 2633  &       & TL    & 2953  &       & TL    & 2903 \\
		\textbf{472} & TL    & 2955  &       & TL    & 3167  &       & TL    & 3194 \\
		\textbf{482} & TL    & 5054  &       & TL    & 3715  &       & TL    & 3885 \\
		\textbf{500} & TL    & 6238  &       & TL    & 5653  &       & TL    & 4278 \\
		\bottomrule
		\bottomrule
	\end{tabular}%
}
}
	\label{tab5}%
\end{table}%

\begin{table}[!h]
	\centering
	\caption{Hyper volume indicator values ($I_H$) for the combination of $\epsilon$-constraint method (e), BB method and the proposed matheuristic (Mat) method with M3 model for $\alpha$-level 0.7 and its corresponding $k$-value for test instances of size 21-500 with the sample size of 10 scenarios.}
	\small\addtolength{\tabcolsep}{-5.5pt}
	\makebox[\linewidth]{
	\scalebox{0.95}{
		
		\begin{threeparttable}
			\begin{tabular}{c|cccccccc}
				\multicolumn{1}{r}{} & \multicolumn{1}{l}{$\alpha$/$k$ value: 0.7/3} &       &       &       &       &       &       &  \\
				\cmidrule{2-9}    \multicolumn{1}{c}{} & \multicolumn{1}{l}{\textbf{e-M3}} &       &       & \multicolumn{1}{l}{\textbf{BB-M3}} &       &       & \multicolumn{1}{l}{\textbf{Mat-M3}} &  \\
				\cmidrule{1-3}\cmidrule{5-6}\cmidrule{8-9}    \multicolumn{1}{c}{\textit{\#Node}} & \textbf{\#NDP} & \textbf{$I_H$} &       & \textbf{\#NDP} & \textbf{$I_H$} &       & \textbf{\#NDP} & \textbf{$I_H$} \\
				\cmidrule{1-3}\cmidrule{5-6}\cmidrule{8-9}    21    & 12    & 3.3095 &       & 12    & 3.3095 &       & 12    & 3.3095 \\
				44    & 29    & 22.8869 &       & 29    & 22.8869 &       & 29    & 22.8869 \\
				56    & 35    & 47.6391 &       & 35    & 47.6391 &       & 35    & 47.6391 \\
				72    & 44    & 72.8961 &       & 44    & 72.8961 &       & 44    & 72.8961 \\
				90    & 53    & 146.6210 &       & 53    & 146.6210 &       & 53    & 146.6200 \\
				106   & 71    & 281.2760 &       & 71    & 281.2760 &       & 71    & 281.1210 \\
				120   & 77    & 306.7360 &       & 77    & 306.7360 &       & 77    & 306.6740 \\
				163   & 101   & 710.9820 &       & 101   & 710.9820 &       & 101   & 710.9820 \\
				182   & 106   & 746.2770 &       & 106   & 746.2770 &       & 106   & 746.2550 \\
				203   & 115   & 896.4520 &       & 115   & 896.4520 &       & 115   & 896.2820 \\
				\textbf{254} & 12    & 117.4060 &       & 120   & 1256.8300 &       & 135   & 1255.8800 \\
				\textbf{275} & 17    & 212.2490 &       & 122   & 1722.6300 &       & 154   & 1723.1900 \\
				\textbf{295} & 2     & 0.0000 &       & 89    & 2440.5200 &       & 155   & 2442.9200 \\
				\textbf{326} & 2     & 0.0000 &       & 35    & 3185.5200 &       & 176   & 3222.0100 \\
				\textbf{355} & 2     & 0.0000 &       & 82    & 4236.5000 &       & 205   & 4250.3400 \\
				\textbf{388} & 2     & 0.0000 &       & 37    & 4335.4900 &       & 218   & 4389.7100 \\
				\textbf{410} & 2     & 0.0000 &       & 80    & 5089.8000 &       & 241   & 5113.6900 \\
				\textbf{436} & 2     & 0.0000 &       & 112   & 5017.4900 &       & 246   & 5025.8700 \\
				\textbf{449} & 2     & 0.0000 &       & 67    & 4838.0800 &       & 255   & 4866.1300 \\
				\textbf{472} & 2     & 0.0000 &       & 39    & 5071.0500 &       & 264   & 5129.3900 \\
				\textbf{482} & 3     & 22.9398 &       & 40    & 5868.9600 &       & 271   & 5960.6000 \\
				\textbf{500} & 2     & 0.0000 &       & 35    & 2591.1700 &       & 181   & 4322.1500 \\
				\bottomrule
				\bottomrule
			\end{tabular}%
			\begin{tablenotes}
				\item Notes. $I_H$ values are in scale of $10^8$.
				\item $^*$ marks the last instance in which the entire exact Pareto frontier is found by e and BB method.
				\item $^+$ indicates the value of $I_H$ where only the reference point (Nadir point) is found.
			\end{tablenotes}
		\end{threeparttable}
	}
}
	\label{tab7}%
\end{table}%

\textcolor{black}{The performance of the generic FPBH method is reported in Table~\ref{tabFPBH}. The results show that the proposed Mat method strongly outperforms the FPBH on our MIP model. The FPBH approach usually reaches the time-limit.}

\begin{table}[!h]
	\arrayrulecolor{black}
	\captionsetup{labelfont={color=black},font={color=black}}
	\centering
	\caption{\textcolor{black}{Performance comparison of the matheuristic (Mat) method and the generic FPBH method with M3 for $\alpha$-level 0.7. T[s] indicates the run time in seconds.}}
	\scalebox{0.85}{
		\begin{threeparttable}
			\begin{tabular}{>{\color{black}}c|>{\color{black}}c>{\color{black}}c>{\color{black}}c>{\color{black}}c>{\color{black}}c>{\color{black}}c>{\color{black}}c}
				\multicolumn{1}{r}{} & \multicolumn{1}{l}{\textcolor{blue}{$\alpha$ value: 0.7}} &       &       &       &       &       &  \\
				\cmidrule{2-8}    \multicolumn{1}{c}{} & \multicolumn{1}{l}{\textcolor{black}{\textbf{Mat-M3}}} &       &       &       & \multicolumn{1}{l}{\textcolor{black}{\textbf{FPBH-M3}}} &       &  \\
				\cmidrule{1-4}\cmidrule{6-8}    \multicolumn{1}{c}{\textcolor{black}{\textit{\#Node}}} & \textbf{T[s]} & {\textbf{\#NDP}} & \textbf{$I_H$} &       & \textbf{T[s]} & \textbf{\#NDP} & \textbf{$I_H$} \\
		\cmidrule{1-4}\cmidrule{6-8}   \textbf{ 21}    & 1     & 12    & 3.3095 &       & 386 & 12    & 3.3095 \\
		\textbf{44}    & 7     & 29    & 39.9519 &       & 7188 & 26    & 39.3993 \\
		\textbf{56}    & 7     & 35    & 84.6040 &       & TL    & 35    & 84.3653 \\
		\textbf{72}    & 25    & 44    & 137.2344 &       & TL    & 48    & 136.9818 \\
		\textbf{90}    & 26    & 53    & 146.5520 &       & TL    & 48    & 133.7343 \\
		\textbf{106}   & 46    & 71    & 281.1210 &       & TL    & 56    & 279.2208 \\
		\textbf{120}   & 45    & 77    & 306.6740 &       & TL    & 52    & 293.1369 \\
		\textbf{163}   & 120   & 101   & 710.9820 &       & TL    & 52    & 704.0039 \\
		\textbf{182}   & 165   & 106   & 746.2550 &       & TL    & 39    & 735.5981 \\
		\textbf{203}   & 356   & 115   & 896.2820 &       & TL    & 42    & 883.1023 \\
		\textbf{254} & 2999  & 135   & 1255.8800 &       & TL    & 17    & 1202.0203 \\
		\textbf{275} & 2612  & 154   & 1723.1900 &       & TL    & 33    & 1619.1607 \\
		\textbf{295} & 1118  & 155   & 2442.9200 &       & TL    & 12    & 2309.4689 \\
		\textbf{326} & 1258  & 176   & 3222.0100 &       & TL    & 6     & 2586.6292 \\
		\textbf{355} & 1487  & 205   & 4250.3400 &       & TL    & 3     & 3144.0494 \\
		\textbf{388} & 2134  & 218   & 4389.7100 &       & TL    & 3     & 1348.2725 \\
		\textbf{410} & 2768  & 241   & 5113.6900 &       & TL    & 2     & 136.4556 \\
		\textbf{436} & 2665  & 246   & 5025.8700 &       & TL    & 3     & 1087.7333 \\
		\textbf{449} & 2992  & 255   & 4866.1300 &       & TL    & 3     & 1075.0511 \\
		\textbf{472} & 3384  & 264   & 5129.3900 &       & TL    & 3     & 1019.8861 \\
		\textbf{482} & 4351  & 271   & 5960.6000 &       & TL    & 2     & 521.7097 \\
		\textbf{500} & 6469  & 181   & 4322.1500 &       & TL    & 9     & 4016.9114 \\
		\bottomrule
	\end{tabular}%
	\begin{tablenotes}
		\item \textcolor{black}{Notes. $I_H$ values are in scale of $10^8$.}
	\end{tablenotes}
\end{threeparttable}
}
\label{tabFPBH}%
\end{table}%

\begin{figure}[!h]
	\arrayrulecolor{black}
	\begin{center}
		\includegraphics[scale=0.8]{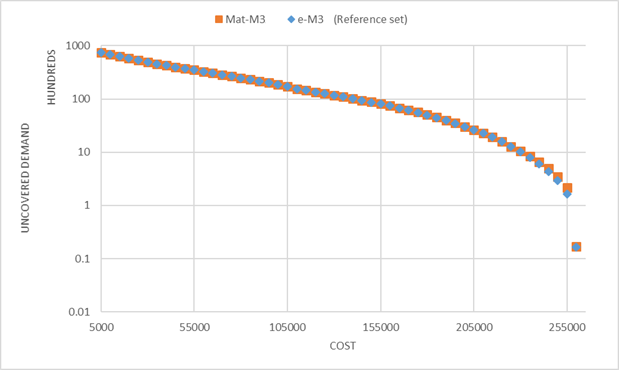}
	\end{center}
	\caption{Logarithmic plot of Pareto frontier for the instance size 90\\
	}
	\label{fig2}
\end{figure}

As an example, the Logarithmic plot of the Pareto frontier resulting from different approaches, for the instance size of 90 nodes which has the worst value of $I_\epsilon$ is shown in Figure~\ref{fig2}. It shows that the Mat method generates almost the same Pareto frontier as the reference set. The main differences are at the right-end of the Pareto frontier, where the location cost is high, and the uncovered demand is at the lower values.

As mentioned before, we compute the stochastic measures in order to demonstrate the effectiveness of incorporating the risk measure of CVaR into the model. We compute the RVPI and RVSS  for different levels of $\alpha$. In Table~\ref{tab11}, we report the relative values of RVPI and RVSS, which are the absolute values divided by the optimal objective value of the RRP (i.e., RVPI/RRP and RVPI/RRP). The results are the average and maximum of the computed values over all non-dominated solutions for the instance size of 21 and 44 nodes with a scenario sample size of 10. According to the results, the value of RVSS is significantly large relative to the RRP value, which indicates that it is worth to solve the risk-averse model with respect to different risk preferences. Moreover, RVPI values are also high, but they get smaller with increasing levels of $\alpha$.

Furthermore, we compute the RVSS values at first by fixing the value of the first-stage decisions ($y_j$) based on the solutions of the deterministic average model (EA) and then, based on the risk-neutral solutions (SP). Table~\ref{tab12} shows the results for the instance size 21 and 44 with a scenario sample size of 10. Comparing the RVSS values indicates that once the randomness has been taken into consideration, the solutions are more robust than the deterministic average model. 

\begin{table}[!h]
	\centering
	\caption{Relative values of stochastic measures (RVPI/RRP, RVSS/RRP) for the intance size 21 and 44 with sample size 10}
	\scalebox{0.85}
	{
	\begin{tabular}{cccccccccccccccc}
		&       & \multicolumn{1}{l}{$\alpha$ value} &       &       &       &       &       &       &       &       &       &       &       &       &  \\
		\cmidrule{3-4}\cmidrule{6-7}\cmidrule{9-10}\cmidrule{12-13}\cmidrule{15-16}          &       & \textbf{0} &       &       & \textbf{0.2} &       &       & \textbf{0.5} &       &       & \textbf{0.7} &       &       & \multicolumn{2}{c}{\textbf{0.9}} \\
		\cmidrule{1-1}\cmidrule{3-16}    \textit{\#Node} &       & \textbf{RVSS} & \textbf{RVPI} &       & \textbf{RVSS} & \textbf{RVPI} &       & \textbf{RVSS} & \textbf{RVPI} &       & \textbf{RVSS} & \textbf{RVPI} &       & \textbf{RVSS} & \textbf{RVPI} \\
		\cmidrule{1-1}\cmidrule{3-16}    \textbf{21} & \textbf{Avg} & 308.31 & 16.86 &       & 303.45 & 14.65 &       & 277.25 & 13.08 &       & 247.55 & 13.08 &       & 264.19 & 7.23 \\
		& \textbf{Max} & 1969.65 & 70.49 &       & 1868.30 & 70.49 &       & 1511.80 & 73.86 &       & 1118.88 & 73.86 &       & 1237.15 & 38.74 \\
		&       &       &       &       &       &       &       &       &       &       &       &       &       &       &  \\
		\textbf{44} & \textbf{Avg} & 646.49 & 22.78 &       & 632.34 & 21.48 &       & 588.77 & 18.14 &       & 501.02 & 15.60 &       & 411.45 & 6.21 \\
		& \textbf{Max} & 7601.92 & 100.00 &       & 7412.21 & 100.00 &       & 6800.00 & 100.00 &       & 5570.52 & 100.00 &       & 3061.54 & 100.00 \\
		\bottomrule
		\bottomrule
	\end{tabular}%
}
	\label{tab11}%
\end{table}%

\begin{table}[!h]
	\centering
	\caption{RVSS with fixed first-stage solutions obtained by solving the risk-neutral model (SP) compared to RVSS with fixed first-stage solutions obtained by solving the deterministic average model (EA) for the intance size 21 and 44 with sample size 10}
	\scalebox{0.85}
	{
	\begin{tabular}{cccccccccccccccc}
		&       & \multicolumn{11}{l}{$\alpha$ value}                                                 &       &       &  \\
		\cmidrule{3-4}\cmidrule{6-7}\cmidrule{9-10}\cmidrule{12-13}\cmidrule{15-16}          &       & \multicolumn{2}{c}{\textbf{0}} &       & \multicolumn{2}{c}{\textbf{0.2}} &       & \multicolumn{2}{c}{\textbf{0.5}} &       & \multicolumn{2}{c}{\textbf{0.7}} &       & \multicolumn{2}{c}{\textbf{0.9}} \\
		\cmidrule{1-1}\cmidrule{3-16}    \textit{\#Node} &       & \textbf{SP} & \textbf{EA} &       & \textbf{SP} & \textbf{EA} &       & \textbf{SP} & \textbf{EA} &       & \textbf{SP} & \textbf{EA} &       & \textbf{SP} & \textbf{EA} \\
		\cmidrule{1-1}\cmidrule{3-16}    \multirow{2}[1]{*}{\textbf{21}} & \textbf{Avg} & 0.00  & 970.47 &       & 4.90  & 1145.94 &       & 30.77 & 1479.12 &       & 63.22 & 1903.77 &       & 88.42 & 2540.69 \\
		& \textbf{Max} & 0.00  & 1247.10 &       & 58.87 & 1483.75 &       & 240.20 & 1877.20 &       & 330.33 & 2402.66 &       & 320.00 & 3293.00 \\
		&       &       &       &       &       &       &       &       &       &       &       &       &       &       &  \\
		\multirow{2}[1]{*}{\textbf{44}} & \textbf{Avg} & 0.00  & 2018.86 &       & 2.82  & 2373.36 &       & 35.98 & 3159.73 &       & 117.69 & 3859.85 &       & 329.62 & 5712.18 \\
		& \textbf{Max} & 0.00  & 2856.50 &       & 42.40 & 3297.00 &       & 173.60 & 4209.00 &       & 440.34 & 5160.97 &       & 868.00 & 7373.00 \\
		\bottomrule
		\bottomrule
	\end{tabular}%
}
	\label{tab12}%
\end{table}%

\begin{figure}[H]
	\begin{center}
		\includegraphics[width=10cm,scale=0.8]{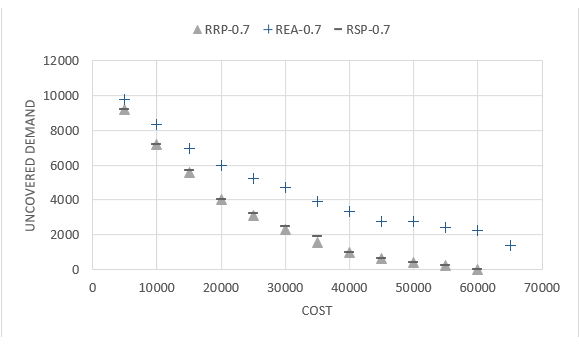}
	\end{center}
	\caption{Pareto frontiers for instance size 21\\
	}
	\label{fig4}
\end{figure}

Figure~\ref{fig4} compares the Pareto frontiers obtained by solving the RRP, the RRP with fixed first-stage solutions from the average model (REA), and the RRP with fixed first-stage solutions from the risk-neutral model (RSP) at $\alpha=0.7$. It shows that incorporating randomness in the model gives a more robust solution. Note that, the REA model obtains solutions at the left-end of the Pareto frontier that are almost as good as RRP and RSP.

\begin{figure}[!h]
	
	\begin{center}
		\includegraphics[width=1\textwidth]{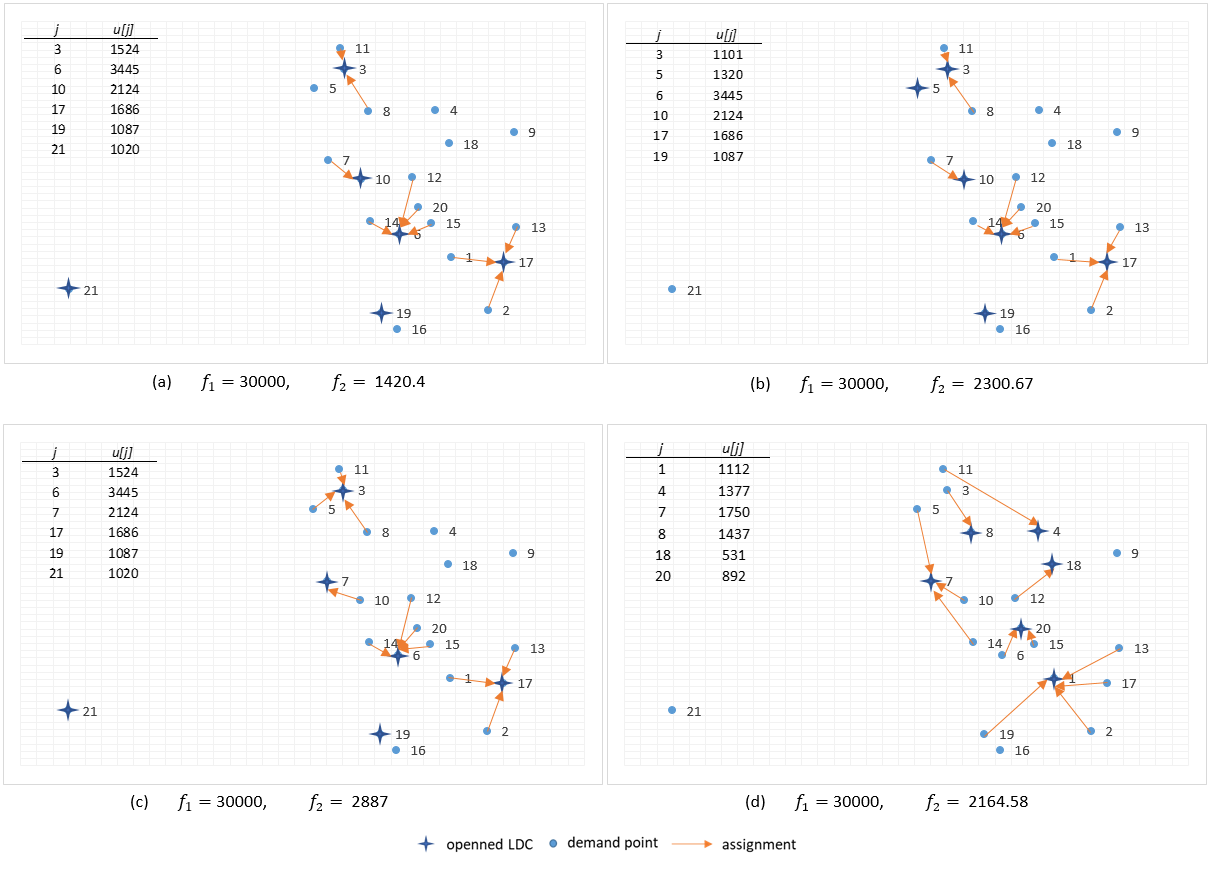}
	\end{center}
	\caption{Comparison of the value decision variables for one the solutions of the efficient set for size 21: (a) $\alpha=0$, (b) $\alpha=0.7$, (c) $\alpha=0.9$ and (d) average deterministic case \\
	}
	\label{fig5}
	
\end{figure}

We also compare the efficient solutions in the solution space for the average deterministic model and the CVaR model with different levels of $\alpha$ for the instance size 21 (Figure~\ref{fig5}). The reported solutions are associated with one point (solution 6 out of 12 Pareto solutions) along the Pareto frontier. A comparison of efficient solutions also highlights the previous finding. It shows that once the randomness in different levels of $\alpha$ is considered, the decisions of the opened PODs ($y_j$) and the allocation of demand points ($x_{ij}$) are more or less similar. Whereas, these decisions in the average deterministic case are entirely different. 

\section{Conclusions}
\label{sec:conclusion}

In this paper, we investigate a two-stage \textcolor{black}{Bi-SSCFL} model \textcolor{black}{extendable to design} a last-mile relief network problem. It incorporates the trade-off between a deterministic objective (minimization of location cost) and an uncertain one (minimization of uncovered demand) where demand values are uncertain.

To find an efficient and reliable methodology to address the problem, we evaluate different methods. We combined two criterion space search frameworks, the $\epsilon$-constraint and the balanced box methods, and three different uncertainty approaches, the widely used two-stage risk-neutral stochastic programming, two-stage adaptive robust optimization and two-stage risk-averse stochastic programming using CVaR. Two variants of linear reformulations of CVaR are taken into consideration.
 
Moreover, we introduce a matheuristic method in order to find high-quality approximations of the Pareto frontier for large-size instances. It is a combination of an adaptive RINS approach, a local branching framework and the $\epsilon$-constraint method.

We perform \textcolor{black}{an in-depth} computational study and compare the results. The experiment illustrates how the proposed matheuristic method not only outperforms the exact frameworks for large-size instances, \textcolor{black}{but also outperforms the FPBH method} \citep{pal2019fpbh}. Additionally, quantifying the VSS measure also shows that incorporating uncertainty into the model gives more robust results in comparison with the deterministic model. In our future research, we would like to extend the model by considering routing constraints where there are multiple sources of uncertainty. There are also potentials for improvement of solution methods by further research on  the combination of MIP and metaheuristics for bi-objective optimization. \textcolor{black}{For example, a decomposition-based approach to deal with the polyhedral subset-based representation of the CVaR might be explored.}

\section*{Acknowledgments}
The authors want to thank \cite{tricoire2012bi} for providing us with the real-world data set of the drought case study. \textcolor{black}{The valuable comments and suggestions of two anonymous reviewers are also greatly appreciated. This research was funded in part by the Austrian Science Fund (FWF) [P 31366-NBL]. For the purpose of open access, the authors have applied a CC BY public copyright license to any Author Accepted Manuscript version arising from this submission.}

\begin{appendices}
	\section{Conditional value-at-risk (CVaR)}
	\label{sec:AppA}
	
	The precise definition of CVaR at confidence level $\alpha \in [0,1) $ of a random variable $X$ is given by \cite{rockafellar2000optimization} as follows:
	
	\begin{align}
		\text{CVaR}_\alpha(X)=\text{min} \{\eta+\frac{1}{1-\alpha}E([X-\eta]_+) : \eta \in \mathbb{R}\},
	\end{align}
	
	Therein, $[b]_+=\text{max}\{b,0\}, b \in \mathbb{R}$ and $\eta$ is the Value-at-Risk, VaR$_\alpha(X)$, of variable $X$ at confidence level $\alpha$. It is defined as follows, where $F_X(.)$ is the cumulative distribution function of a random variable $X$. 
	
	\begin{align}
		\text{VaR}_\alpha(X)=\text{min} \{\eta \in \mathbb{R}: F_X(\eta)\geq \alpha\},
	\end{align}
	
	CVaR$_\alpha(X)$ can be interpreted as the conditional expected value exceeding the VaR at the confidence level $\alpha$ (see, \cite{rockafellar2000optimization}).
	
	The auxiliary real variables $w_s$ for $s=1,\dots,N$ are introduced to reduce (28) to a linear programming problem:
	
	\begin{align}
		&\text{min}& &\text{CVaR}_\alpha(X)=\eta+\frac{1}{1-\alpha}\frac{1}{N}\sum\limits_{s=1}^N  w_s&\\
		&\text{s.t.}& &w_s\geq X-\eta& &\forall s=1,\dots,N&\\
		& & &w_s\geq 0& &\forall s=1,\dots,N&\\
		& & &\eta \in \mathbb{R}& &&
	\end{align}

	\section{Multi-objective concepts and definitions}
	\label{sec:AppB}
	
	Due to the bi-objective nature of our problem, it is not possible to find a unique optimal solution, but rather the so-called set of efficient solutions. A MOO problem is given as follows:
	
	\begin{align}
		&\min_{x \in X}& &f(x)=(f_1(x),f_2(x),...,f_m(x))&
	\end{align}
	
	Therein, $m \geq 2$ is the number of objectives, $x\in X$, where $X$ is the set of feasible solutions in the \textit{decision space} and functions $f_i: X \to \mathbb{R}$ are the objective functions, and $Z = f(X)$ the image of the set of feasible solutions in the \textit{criterion space}.
	
	\begin{definition}{Pareto dominance}\\
		Solution $x^* \in X$ dominates $x \in X$ if and only if $x^*$ is as good as $x$ with respect to all objectives, and better than $x$ with respect to at least one of the objectives, i.e.,
		\begin{align*}
			\begin{cases}
				f_i(x^*) \leq f_i(x) \quad \text{for all   } i\in \{1,...,m\}\\
				f_j(x^*) < f_j(x) \quad \text{for at least one  } j\in \{1,...,m\}
			\end{cases}\\
		\end{align*}
	\end{definition}
	
	\begin{definition}{Pareto frontier}\\
		$x \in X$ is an efficient solution of MOO (34) if there is no $x^* \in X$ that dominates $x$. The set of efficient solutions of MOO is denoted by $X_e$. The objective vector $f(x)$ of an efficient solution $x \in X_e$ is called non-dominated solution, and $Z=f(X_e)=\{f(x): x \in X_e\}$ is the set of non-dominated objective points (NDP) or Pareto frontier.
		
	\end{definition}

\section{Results for larger number of scenarios}
\label{sec:AppC}

\begin{table}[H]
	\centering
	\caption{Performance comparison of a combination of $\epsilon$-constraint method (e), and the proposed matheuristic (Mat) method with M3 for $\alpha$-level 0.7 and its corresponding $k$-value for larger number of scenarios. T[s] indicates the run time in seconds.}
	\small\addtolength{\tabcolsep}{-5.5pt}
	\makebox[\linewidth]{
		\scalebox{0.95}{
			
			\begin{threeparttable}
				
					\begin{tabular}{ccccccccccc}
						& \multicolumn{10}{l}{\textbf{0.7/k}} \\
						\cmidrule{2-11}          &       &       &       & \multicolumn{3}{c}{\textbf{e-M3}} &       & \multicolumn{3}{c}{\textbf{Mat-M3}} \\
						\cmidrule{1-7}\cmidrule{9-11}    \textit{\#Node} & \textit{\#Scenario} & \textit{k} &       & T[s]  & \#NDP & $I_H$   &       & T[s]  & \#NDP & $I_H$ \\
						\cmidrule{1-7}\cmidrule{9-11}    \multirow{3}[2]{*}{21} & 100   & 30    &       & 161   & 17    & 6.0735 &       & 27    & 17    & 6.0527 \\
						& 500   & 150   &       & 6476  & 16    & 5.5603 &       & 348   & 16    & 5.5600 \\
						& 1000  & 300   &       & TL    & 6     & 2.0746 &       & 1668  & 17    & 6.0657 \\
						\cmidrule{1-7}\cmidrule{9-11}    \multirow{3}[2]{*}{44} & 100   & 30    &       & 1076  & 30    & 25.3530 &       & 178   & 30    & 25.3530 \\
						& 500   & 150   &       & TL    & 5     & 3.5413 &       & 1374  & 32    & 27.4957 \\
						& 1000  & 300   &       & TL    & 3     & 1.2212 &       & TL    & 3     & 1.2212 \\
						\cmidrule{1-7}\cmidrule{9-11}    \multirow{3}[2]{*}{56} & 100   & 30    &       & 2267  & 39    & 55.3764 &       & 234   & 39    & 55.3543 \\
						& 500   & 150   &       & TL    & 3     & 2.0392 &       & 2675  & 40    & 56.5806 \\
						& 1000  & 300   &       & TL    & 3     & 1.9092 &       & TL    & 2     & 0.0000 \\
						\cmidrule{1-3}\cmidrule{5-7}\cmidrule{9-11}    \multirow{3}[2]{*}{72} & 100   & 30    &       & 4186  & 50    & 101.2260 &       & 180   & 50    & 101.2260 \\
						& 500   & 150   &       & TL    & 3     & 2.7965 &       & TL    & 2     & 0.0000 \\
						& 1000  & 300   &       & TL    & 2     & 0.0000 &       & TL    & 2     & 0.0000 \\
						\bottomrule
						\bottomrule
					\end{tabular}%
				
				\begin{tablenotes}
					\item \textcolor{black}{Notes. $I_H$ values are in scale of $10^8$.}
				\end{tablenotes}
			\end{threeparttable}
		}
	}
	\label{tab6}%
\end{table}%
	
\end{appendices}

\linespread{1}
\small 

\bibliography{paperRef}


\end{document}